\documentclass[10pt]{amsproc}
\usepackage{amsmath,amsthm,amsfonts,amscd,float,amssymb}
\newcommand{\Gal}{{\rm Gal}}

\newcommand{\CC}{\mathbb{C}}
\newcommand{\KK}{\mathbb{K}}
\newcommand{\PP}{\mathbb{P}}
\newcommand{\FF}{\mathbb{F}}
\newcommand{\OO}{\mathcal{O}}
\newcommand{\ZZ}{\mathbb{Z}}
\newcommand{\QQ}{\mathbb{Q}}

\newcommand{\bx}{\mathbf{x}}

\newcommand{\pp}{\mathfrak{p}}

\newcommand{\ba}{\mathbf{a}}
\newcommand{\RR}{\Bbb{R}}
\newcommand{\Frob}{{\rm Frob}\thinspace}
\newcommand{\Norm}{{\rm Norm}\thinspace}

\newcommand{\trace}{{\rm tr}\thinspace}

\newcommand{\Ind}{{\rm Ind}\thinspace}
\newcommand{\Res}{{\rm Res}\thinspace}

\theoremstyle{plain}
\newtheorem{thm}{Theorem}[section]

\theoremstyle{definition}
\newtheorem{defn}{Definition}[section]
\newtheorem{conj}{Conjecture}[section]

\newtheorem{question}{Question}[section]

\theoremstyle{definition}

\floatstyle{plain}
\newfloat{Program}{thp}{lop}[section]
\newfloat{Table}{thp}{lop}[section]

\begin{document}

\title{Update on Modular Non-Rigid Calabi-Yau Threefolds}
\author{Edward Lee}
\address{Department of Mathematics, UCLA, Los Angeles, CA 90095}
\email{lee@math.ucla.edu}
\begin{abstract}
We review some recent results on the modularity of non-rigid Calabi-Yau threefolds.
\end{abstract}
\maketitle


\section{Introduction}
\subsection{Calabi-Yau varieties}

This note can be considered a sequel to a previous survey article by Yui \cite{bib:Yui}.  The main objects of discussion are Calabi-Yau varieties and their cohomological $L$-series.  We will review recent constructions of Calabi-Yau varieties which are proven to be modular in the sense that their cohomological $L$-series are associated to modular forms.

We recall first some basic definitions.

\begin{defn}
Let $X$ be a complex projective variety of dimension $d$.  $X$ is a {\em Calabi-Yau variety} if it satisfies the conditions
\end{defn}

\begin{enumerate}

\item $H^i(X, \OO_X) = 0$ for $0 < i < d$;

\item The canonical sheaf $K_X$ of $X$ is trivial.

\end{enumerate}

The Hodge cohomology groups and Hodge numbers are defined by

\begin{equation*}
\begin{aligned}
H^{i,j}(X) &= H^j(X, \Omega^i_X), h^{i,j}(X) = \dim H^{i,j}(X).
\end{aligned}
\end{equation*}

If $X$ is smooth, the cohomology groups of $X$ admit the Hodge decomposition

\begin{equation*}
\begin{aligned}
H^{k}(X,\CC) = \oplus_{i+j = k} H^{i,j}(X).
\end{aligned}
\end{equation*}

Furthermore, by complex conjugation we have $\overline{H^{i,j}(X)} \cong H^{j,i}(X)$, and by Serre duality and the triviality of the canonical sheaf we have $H^{d-i,d-j}(X) \cong H^{i,j}(X)$.

The Betti numbers $b^{i}(X)$ and the Euler characteristic $\chi(X)$ are defined by $b^{i}(X) = \dim H^i(X, \CC)$ and

\begin{equation*}
\begin{aligned}
\chi(X) &= \sum_{i=0}^{2d} (-1)^i b^i(X). \\
\end{aligned}
\end{equation*}

Smooth one-dimensional Calabi-Yau varieties are simply elliptic curves, with Hodge diamond

\begin{equation*}
\begin{tabular}{ccc}
         & 1 &  \\
        1 &  & 1 \\
         & 1 &  \\
\end{tabular}
\end{equation*}

Smooth two-dimensional Calabi-Yau varieties are $K3$ surfaces, with Hodge diamond

\begin{equation*}
\begin{tabular}{ccccc}
       &   & 1 &   & \\
       & 0 &   & 0 & \\
     1 &   & 20&   &  1 \\
       & 0 &   & 0 & \\
       &   & 1 &   & \\
\end{tabular}
\end{equation*}

In this article, we will be most interested in three-dimensional Calabi-Yau varieties.  Smooth Calabi-Yau threefolds have Hodge diamond

\begin{equation*}
\begin{tabular}{ccccccc}
    &  &   & 1 &   & & \\
    &  & 0 &   & 0  & & \\
   & 0 &  & $h^{11}$ &  & 0 & \\
  1&  & $h^{12}$  & & $h^{12}$  &   & 1\\
   & 0 &  & $h^{11}$  &  & 0 & \\
     & & 0  &  & 0  &  & \\
     & &   &  1 &   &  & \\
\end{tabular}
\end{equation*}

though we shall have occasion to discuss results in other dimensions.

Calabi-Yau threefolds can be classified as {\em rigid} or {\em non-rigid} according to whether $h^{2,1}$ is zero or nonzero respectively.  By Kodaira-Spencer theory, the infinitesimal deformations of $X$ are classified by $H^1(X, T_X)$, which is isomorphic to $H^1(X, \Omega^2_X)$ by the triviality of $K_X$.  A rigid Calabi-Yau threefold has $h^{2,1} = 0$ and thus admits no infinitesimal deformations, whence the name.  Since the infinitesimal deformations of any Calabi-Yau variety $X$ are unobstructed \cite{bib:Tian} \cite{bib:Todorov}, $X$ moves in a $h^{2,1}$-dimensional family.

\subsection{$L$-series of a Calabi-Yau variety}

We now introduce the cohomological $L$-series of a Calabi-Yau variety $X$, which encodes a vast amount of arithmetic information about $X$.

Suppose that the Calabi-Yau variety $X$ is defined over a number field $K$, and let $\OO_K$ be the ring of integers of $K$.  Then $X$ always has a model defined over $\OO_K$.  Fix some model over $\OO_K$, and let $\pp$ be a prime ideal of $\OO_K$.  We say that the prime $\pp$ is a {\em good} prime of $X$ if the reduction $X_\pp$ is smooth; otherwise we say that it is a {\em bad} prime of $X$.  Let $q$ be the norm of the prime ideal $\pp$ and let $\FF_q$ be the field of $q$ elements.

If $\pp$ is a good prime of $X$, let $\Frob_{\pp}$ denote the Frobenius morphism on $X_{\pp}$ induced by the coordinate map $x \rightarrow x^q$.  Then $\Frob_{\pp}$ acts on the \'etale cohomology groups $H^i_{et}(\overline{X}_{\pp}, \QQ_l)$, where $\overline{X}_{\pp}$ denotes base change of $X_{\pp}$ to $\overline{\FF}_q$ and $l$ is a rational prime which is coprime to $q$.  By abuse of notation we will often continue to write $X_{\pp}$ for its base change to $\overline{\FF}_\pp$.

The Frobenius morphism induces an endomorphism of the \'etale cohomology groups, which we denote by $\Frob^{*}_{\pp}$.  Let

\begin{equation*}
\begin{aligned}
P_i^{\pp}(T) &= \det(1 - \Frob^{*}_{\pp} T | H^i_{et}(\overline{X}, \QQ_l)) \\
\end{aligned}
\end{equation*}

be the characteristic polynomial of the Frobenius endomorphism on $H^i_{et}(\overline{X}, \QQ_l)$.  We will be most interested in the middle-dimensional cohomology group of $X$.

If $\pp$ is a prime of bad reduction for $X$, then we can still define the polynomial

\begin{equation*}
\begin{aligned}
P_i^{\pp}(T) &= \det (1 - \Frob^{*}_{\pp} T | H^i_{et}(\overline{X}, \QQ_l)^{I_{\pp}}) \\
\end{aligned}
\end{equation*}

where $I_{\pp}$ is the inertia group at $\pp$.  However, unlike the elliptic curve case, our current understanding of bad reduction of higher dimensional Calabi-Yau varieties is very thin.  In most of the examples below we will determine $L$-series modulo the Euler factors at bad primes.

\'Etale cohomology is a Weil cohomology theory, and thus satisfies a Lefschetz fixed-point theorem:

\begin{thm}
If $f$ is an endomorphism of the smooth variety $X$, then the number of fixed points of $f$ satisfies the formula

\begin{equation}
Fix(f,X) = \sum_{i=0}^{2d} (-1)^i \trace f^{*}(H^i_{et}(X,\QQ_l))
\end{equation}
\end{thm}

In the case $f = \Frob_\pp$, the number of fixed points will simply be the number of $\FF_q$-rational points of $X$, and this equals the alternating sum of the traces of the Frobenius action on the cohomology groups.  In the Calabi-Yau threefold case, the traces of $\Frob_\pp$ on $H^0$ and $H^6$ are equal to 1 and $q^3$ respectively, and $H^1 = H^5 = 0$.  If one understands the Frobenius action on the groups $H^2$ and $H^4$, then the trace of Frobenius on $H^3$ can be determined by counting points (often using a computer).  For example, if $H^2$ is spanned by algebraic cycles defined over $\OO_K/\pp$, then Frobenius acts by multiplication by $q$ on the whole space of algebraic cycles, and thus the trace of Frobenius on $H^2$ is just $qb^2(X)$.

Furthermore, Deligne's proof of the Weil conjectures tells us that

\begin{thm}
$P_i^{\pp}(T)$ is a polynomial with integer coefficients, whose reciprocal roots are all algebraic integers with absolute value $p^{\frac{i}{2}}$.
\end{thm}

The absolute Galois group $G(\overline{\KK}/\KK)$ acts on $X$.  If $l$ is a prime, we have the $l$-adic Galois representation

\begin{equation*}
\rho_{X,l}^{i}: Gal(\overline{\KK}/\KK) \rightarrow GL(H^i_{et}(X, \QQ_l))
\end{equation*}

and we define its $L$-series $L(\rho_{X,l}^{i}, s) = L(H^{i}_{et}(X,\QQ_l))$:

\begin{defn}
Let $X$ be a Calabi-Yau variety of dimension $d$ defined over a number field $\KK$.  For each $i$, the $i$th cohomological $L$-series of $X$ is defined by the Euler product

\begin{equation*}
\begin{aligned}
L(H^{i}_{et}(X, \QQ_l),s) &= \prod_{\pp} P_i^{\pp}(q^{-s})^{-1}, q = \Norm \pp
\end{aligned}
\end{equation*}

\end{defn}

However, in the sequel the factors corresponding to primes of bad reduction $\pp$ for $X$ will often be omitted.

\begin{defn}
A Calabi-Yau variety $X$ of dimension $d$ defined over $\QQ$ is {\em modular} if there is a modular form (or more generally an automorphic form) such that its associated Galois representation is equivalent to $\rho_{X,l}^{d}$.
\end{defn}

In particular, the $L$-series of $X$ coincides with the $L$-series of the modular or automorphic form, possibly up to a finite number of Euler factors.

\section{Elliptic curves}

The dimension 1 case is of course the case of elliptic curves, and the question of modularity of elliptic curves over $\QQ$ was famously settled by Wiles and his collaborators \cite{bib:BCDT} \cite{bib:TW}.  Here the only nontrivial cohomology group is $H^1_{et}(E, \QQ_l)$, and the Lefschetz formula gives us

\begin{equation*}
a_p = \trace (Frob_p^{*}|H^1_{et}(E_p, \QQ_l)) = p + 1 - \# E_p.
\end{equation*}

The conductor $N$ of $E$ is divisible exactly by the bad primes of $E$ and measures the extent to which $E_p$ is singular:

\begin{equation*}
N = \prod_p p^{f_p}
\end{equation*}

where

\begin{equation*}
f_p = 
\begin{cases}
0 & \text{if E has good reduction at $p$}, \\
1 & \text{if $E_p$ has a node}, \\
2 & \text{if $E_p$ has a cusp and $p \neq 2,3$,} \\
2 + \delta_p & \text{if $E_p$ has a cusp and $p = 2$ or $p = 3$} \\
\end{cases}
\end{equation*}

and $\delta_p$ measures the wild ramification in the action of the inertia subgroup $I_p$ of $\Gal(\overline{\QQ}/ \QQ)$ on $H^1_{et}(E,\QQ_l)$.

The modularity theorem then states that there exists a weight 2 cusp form $f$ of level $N$ such that

\begin{equation*}
L(E,s) = L(f,s).
\end{equation*}

\section{K3 surfaces}

In dimension 2, the Calabi-Yau varieties are exactly the $K3$ surfaces.  The only nontrivial cohomology group is $H^2$, which has dimension 22.

Let $NS(X)$ be the Neron-Severi group of $X$, i.e. the group of divisors on $X$ modulo algebraic equivalence; its rank is called the {\em Picard number} of $X$.  Since the divisors of $X$ live in $H^{11}(X) \cap H^2(X,\ZZ)$, the maximum Picard number of $X$ is 20.  $K3$ surfaces achieving this maximum are called {\em singular} $K3$ surfaces.

If $X$ is a singular $K3$ surface, the lattice $H^2(X,\ZZ)$ under cup product decomposes as $NS(X) \oplus T(X)$ (up to finite index), where $T$ is called the {\em transcendental lattice} of $H^2(X,\ZZ)$.  This induces a decomposition of the $L$-series of $X$ as

\begin{equation*}
L(X,s) = L(NS(X) \otimes \QQ_l, s) L(T(X) \otimes \QQ_l,s).
\end{equation*}

\begin{thm}
Let $X$ be an singular $K3$ surface over $\QQ$.  Up to Euler factors at the bad primes, the $L$-series of $T(X)$ is given by

\begin{equation*}
L(T(X) \otimes \QQ_l,s) = L(f,s)
\end{equation*}

where $f$ is a cusp form of weight 3 on a congruence subgroup of $PSL_2(\ZZ)$.  The level $N$ can be computed explicitly.  The algebraic part $L(NS(X),s)$ depends on the fields of definition of the divisors of $X$.
\end{thm}

The proof of this result depends on the work of Inose and Shioda on the classification of singular $K3$ surfaces \cite{bib:IS} and on the work of Livn\'e on the modularity of motivic orthogonal two-dimensional Galois representations \cite{bib:Livne1}.

If the Picard number of $X$ is less than 20, much less is known about the modularity of $X$.  In this case the transcendental part $T(X)$ of $H^2$ has dimension higher than 2.

\section{Modularity of rigid Calabi-Yau threefolds}

Perhaps the first example of a rigid Calabi-Yau threefold which was shown to be modular is a resolution $\tilde{Q}$ of the Schoen quintic threefold $Q = \{ \bx \in \PP^4 : x_0^5 + x_1^5 + x_2^5 + x_3^5 + x_4^5 - 5x_0x_1x_2x_3x_4 = 0 \}$; in \cite{bib:Schoen1} it was shown that the middle cohomology of $\tilde{Q}$ is associated to the weight 4 form 

\begin{equation*}
\eta(z)^4 \eta(5z)^4 + 5 \eta(z)^3 \eta(5z)^4 \eta(25z) + 20 \eta(z)^2 \eta(5z)^4 \eta(25z)^2 + \eta(z) \eta(5z)^4 \eta(25z)^3 + \eta(5z)^4 \eta(25z)^4
\end{equation*}

of level 25, where $\eta$ is the Dedekind eta-function.  Since then, many clever constructions of other rigid modular Calabi-Yau threefolds have been found; for reasons of space we must refer the reader to \cite{bib:Meyer}.

In the rigid case, it is expected that $H^3$ will always correspond to a weight 4 modular form.  Under mild conditions, Dieulefait and Manoharmayum have proven this to be true:

\begin{thm}
\cite{bib:Dieu} \cite{bib:DM}
Suppose that $X$ is a rigid Calabi-Yau threefold over $\QQ$.  Assume that one of the following statements is true:

\begin{enumerate}

\item $X$ has good reduction at 3 or 7 or

\item $X$ has good reduction at 5  or

\item $X$ has good reduction at 3 and the trace of $\Frob_3$ on $U$ is not divisible by 3.

\end{enumerate}

Then $X$ is modular.

\end{thm}

We might hope that the conditions of good reduction at the specific primes 3,5,7 can eventually be removed, as in the gradual proof of modularity for elliptic curves.  

\section{Modular threefolds containing elliptic ruled surfaces}

We now turn our attention to modular non-rigid Calabi-Yau threefolds; the first such example discovered can be found in \cite{bib:NvG}.  

The simplest cases in which the modularity of a nonrigid Calabi-Yau threefold $X$ can be proven occur when the semisimplification of $H^3(X)$ decomposes as $U \oplus V$.  Here $U$ is two-dimensional and $V$ is of the form $\oplus_{i = 1}^{h^{12}} H^3(S_i)$, where the $S_i$ are ruled surfaces over elliptic curves.  (The surfaces $S_i$ may be said to {\em span} $H^{21}(X) \oplus H^{12}(X)$ in this case).  If the surfaces $S_i$ are defined over $\QQ$ (and even if only the scheme-theoretic union $\cup S_i$ is defined over $\QQ$), then it essentially follows from the theorem of Dieulefait-Manoharmayum that $U$ is modular:

\begin{thm}
\cite{bib:HV2}  Let $X$ be a smooth projective threefold defined over $\QQ$ such that $h^{30}(X) = h^{03}(X) = 1$ and assume that $X$ contains birational ruled surfaces $S_i, i = 1, \dots, h^{12}(X)$ which are defined over $\QQ$ and which span $H^{21}(X) \oplus H^{12}(X)$.  Let $\rho$ be the 2-dimensional Galois representation given by $U$, and suppose we have an exact sequence

\begin{equation*}
0 \rightarrow U \rightarrow H^3_{et}(X, \QQ_l) \rightarrow \oplus_{i=1}^{h^{12}} H^3_{et}(S_i, \QQ_l) \rightarrow 0.
\end{equation*}

Assume that one of the following conditions holds:

\begin{enumerate}

\item $X$ has good reduction at 3 or 7 or

\item $X$ has good reduction at 5 or

\item $X$ has good reduction at 3 and the trace of $\Frob_3$ on $U$ is not divisible by 3.

\end{enumerate}

Then $X$ is modular, and the $L$-function of $X$ is given by

\begin{equation}
L(X,s) \circeq L(f_4,s) \prod_{i=1}^{h^{12}} L(g_2^i, s-1)
\end{equation}

where $f_4$ is a weight 4 modular form, the $g_2^i$ are the weight 2 forms associated to the base curves $E_i$, and $\circeq$ denotes equality up to finitely many Euler factors.
\end{thm}

Examples of modular threefolds such that $H^3(X)$ decomposed in this manner were constructed by Hulek-Verrill, Cynk-Meyer, Sch\"utt and Lee.  In each of these cases $U$ is a two-dimensional Galois representation whose modular form can be determined by counting points on $X$ over finite fields.

\subsection{Hypersurfaces in toric varieties}

In \cite{bib:HV1} Hulek and Verrill consider a family $\{ X_{\bf a}, {\bf a} \in \PP^5 \}$ of Calabi-Yau threefolds which are obtained from anticanonical hypersurfaces in toric varieties; they consider a reflexive polytope obtained from the $A_4$ root lattice.  

Let $M_{A_4}$ be the sublattice of points $(x_1,x_2,x_3,x_4,x_5) \in \ZZ^5$ such that $\sum x_i = 0$, and define the polytope $\Delta_{A_4}$ to be the convex hull in $M \otimes \RR$ of the points $(e_i - e_j), 1 \leq i \neq j \leq 5$.  They then take the dual polytope $\Delta^{*}$, consider a subdivision $\tilde{\Delta}^{*}$ and consider hypersurfaces in the resulting toric variety $\tilde{P}$.  The toric variety $\tilde{P}$ contains an open 4-torus $T$, which we describe as $\{ {\bf X} \in \PP^4 : \prod X_i \neq 0 \}$.

Let $X_{\ba}$ be the anticanonical hypersurface of $\tilde{P}$ whose equation, when restricted to the open torus $T$, is

\begin{equation}
(X_1 + \dots + X_5)(\frac{a_1}{X_1} + \dots + \frac{a_5}{X_5}) = t
\end{equation}

for $\ba = (a_1 : a_2 : a_3 : a_4 : a_5 : t) \in \PP^5$.  In general, $X_{\ba}$ has 30 ordinary double point singularities, and special values of $\ba$ yield threefolds with extra nodes, so one takes a big resolution $\tilde{X}_{\ba}$.  (Big resolutions of ordinary double point singularities are non-crepant resolutions, so $\tilde{X}_{\ba}$ is only birationally Calabi-Yau.)

For $\ba$ = $(1:1:1:1:1:25)$,$(1:1:1:9:9:9)$ and $(1:1:4:4:4:16)$, $\tilde{X}_{\ba}$ is shown to contain ruled surfaces $S_i$ over elliptic curves $E_i$ as subvarieties; these elliptic curves $E_i$ are defined over $\QQ$ and their modular forms can be determined explicitly.

The inclusion map $i: \cup S_i \hookrightarrow \tilde{X}_{\bf a}$ is proven to induce an embedding $i_{*}: \oplus_i H_3(S_i) \hookrightarrow H_3(\tilde{X}_{\bf a})$.  By using Poincar\'e duality and using comparison theorems to pass to \'etale cohomology, we get the desired splitting of $H^3(\overline{X_{\ba}})$ (up to semisimplification)

\begin{equation*}
H^3_{et}(X, \QQ_l) \cong \oplus H^3_{et}(S_i, \QQ_l) \oplus U
\end{equation*}

It remains to prove that $U$ is modular.  Now the trace of $\Frob_p$ on $U$ will simply be the trace of $\Frob_p$ on $H^3_{et}(X)$ minus the trace on the $H^3_{et}(S_i)$.  The traces on $H^3_{et}(S_i)$ will just be $p$ times the coefficients of the modular forms of the elliptic curves $E_i$, and the traces of $\Frob_p$ on $\tilde{X}$ are obtained by counting points and using the Lefschetz theorem as described earlier.  Hence the traces of Frobenius on $U$ can be effectively computed.  By inspection, they find a candidate modular form of weight 4 whose coefficients agree with the traces of Frobenius for many values of $p$.  They then use a well-known theorem of Livn\'e (following ideas of Faltings and Serre) which essentially states that two two-dimensional Galois representations $U_1$ and $U_2$ with even trace can be proven to be isomorphic if their traces are the same at some finite set of primes $T$, where $T$ depends on the bad primes of $U_1$ and $U_2$:

\begin{thm}
\cite{bib:Livne1}
Let $\rho_1$ and $\rho_2$ be two 2-adic 2-dimensional Galois representations, unramified outside a finite set of primes $S$.  Let $K_S$ be the compositum of all quadratic extensions of $\QQ$ ramified at primes in $S$, and let $T$ be a set of primes disjoint from $S$.  Then if
\begin{enumerate}

\item$ \trace \rho_1 \equiv \trace \rho_2 \equiv 0 \: and \det \rho_1 \equiv \det \rho_2 ({\rm mod} 2)$;

\item $\{ \Frob_p|{K_S} : p \in T \}$ is equal to the set $\Gal(K_S/K) - Id$;

\item for all $p \in T$, $\trace \rho_1(\Frob_p) = \trace \rho_2(\Frob_p)$;

then $\rho_1$ and $\rho_2$ have isomorphic semisimplifications.

\end{enumerate}

\end{thm}

The upshot of this theorem is that checking that $U$ and the Galois representation associated to a modular form $g$ are isomorphic reduces to the problem of checking that their coefficients are even and that the traces of $U$ are equal to the coefficients of $g$ at some {\em finite} set of primes $T$.

Thus we have

\begin{thm}
For the special values of $\ba$ listed above, the threefold $\tilde{X}_{\ba}$ is modular, with $L$-series equal to

\begin{equation}
L(H^3_{et}(\tilde{X}_{\ba}), s) \circeq L(g,s) \prod_{i=1}^{h^{12}} L(f_i,s-1)
\end{equation}

where the $f_i$ are modular forms of elliptic curves over $\QQ$ and $g$ is some modular form of weight 4.

\end{thm}

\subsection{Products of elliptic fibrations}

In \cite{bib:HV2} and \cite{bib:Schutt} Hulek-Verrill and Sch\"utt consider fiber products of relatively minimal, regular, rational elliptic surfaces with section over $\PP^1$ (in fact, the previous example is birational to a threefold of this type).  Given two such surfaces $\pi_1: Y_1 \rightarrow \PP^1$ and $\pi_2: Y_2 \rightarrow \PP^1$, it was first shown in \cite{bib:Schoen2} that if the fibers of both $Y_1$ and $Y_2$ over the set of common cusps are all of type $I_n$, the fiber product $W = Y_1 \times_{\PP^1} Y_2$ of two such elliptic fibrations is a nodal Calabi-Yau threefold.

Assume now that all degenerate fibers of $\pi_1$ and $\pi_2$ are of type $I_n$, and consider a point $p \in \PP^1$ which is a cusp for exactly one of the fibrations, say $\pi_1$.  Then the fiber $W_p$ is a product of an elliptic curve $E$ and a degenerate fiber of type $I_n$, i.e. a union of elliptic ruled surfaces.  It can be shown by intersection theory that the inclusion induces an embedding

\begin{equation*}
\i_{*}: \oplus_{i=1}^{h^{12}} H_3(E_i \times \PP^1) \rightarrow H_3(X).
\end{equation*}

After using duality and comparison theorems, we obtain a surjection

\begin{equation*}
\i^{*}:   H^3_{et}(X, \QQ_l) \rightarrow \oplus_{i=1}^{h^{12}} H^3_{et}(E_i \times \PP^1, \QQ_l)
\end{equation*}

whose kernel is 2-dimensional, so we are in the situation of Theorem 5.1.

In the examples of \cite{bib:HV2} the representation $U$ has even trace, so once a candidate modular form is found, the theorem of Faltings-Serre-Livn\'e can be used to establish modularity.  

In the examples of \cite{bib:Schutt}, the representations $U$ do not have even trace, so the argument to prove the equivalence of $U$ with the Galois representation of a modular forms $f$ is more subtle.  Here Sch\"utt works out in detail a method first suggested by Serre in \cite{bib:Schoen1} to prove that two two-dimensional 2-adic representations $\rho_1$ and $\rho_2$ are isomorphic.  The idea is to first prove that the mod 2 reductions are isomorphic and absolutely irreducible.  Then under the assumption (to be contradicted) that $\rho_1$ and $\rho_2$ are isomorphic mod $2^\alpha$ but not mod $2^{\alpha+1}$, one constructs an auxiliary Galois representation $\tilde{\rho}: \Gal(\overline{\QQ}/\QQ) \rightarrow M_2(\FF_2) \rtimes GL_2(\FF_2)$ which measures the degree to which $\rho_1$ and $\rho_2$ would not be isomorphic.  Carrying out this method requires one to know the Galois extensions $K$ of $\QQ$ with Galois group $S_3$ or $C_3$ unramified outside $S$ (one way is to check the tables in \cite{bib:Jones}), then to know the Galois extensions $L/K$ with $Gal(L/\QQ)$ isomorphic to $S_3 \times C_2$ or $S_4$ and unramified outside $S$.

A particularly interesting example in \cite{bib:Schutt} is of a Calabi-Yau threefold $\tilde{W_2}$ which has bad reduction at 2,5,11 and 359.  Numerical evidence suggests that the non-elliptic part of its third cohomology should be associated with a newform $f_{55}$ of weight 4 and level 55, but unfortunately this cannot be proven yet due to the difficulty of computing all cubic and quartic Galois extensions of $\QQ$ unramified outside $\{2,5,11,359 \}$.  However, assuming that modularity can be proven, this would be an example where the level of the modular form is not divisible by all the primes of bad reduction, unlike the case of elliptic curves.

\subsection{A Horrocks-Mumford quintic}

The Horrocks-Mumford vector bundle $HM$ is essentially the only known rank 2 indecomposable vector bundle on $\PP^4$, introduced in \cite{bib:HM} and studied extensively since then (see \cite{bib:Hulek} for an introduction).  It was shown in \cite{bib:HM} that the zero locus of an arbitrary section of $HM$ is an abelian surface of degree 10.  Now $\wedge^2 HM \cong \OO(5)$, so a determinantal section of the bundle $\wedge^2 HM$ is a (singular) quintic threefold which is the total space of a pencil of abelian surfaces.

In \cite{bib:Lee} the Horrocks-Mumford quintic $X$ with equation

\begin{equation*}
\sum_{j \in \ZZ/5} (z_j^3 z_{j+1} z_{j-1} - z_j^3 z_{j+2} z_{j-3} - z_j z_{j+1}^2 z_{j-1}^2 + z_j z_{j+2}^2 z_{j-2}^2) = 0
\end{equation*}

was studied (where the indices are read modulo 5).  It was shown that the singular locus of $X$ contains the elliptic normal curves $E_1$ and $E_2$, where $E_1$ is cut out by the five equations $i z_j^2 + z_{j+1}z_{j-1} + z_{j+2}z_{j-2} = 0$ and $E_2$ is the complex conjugate of $E_1$, and there exists a big resolution $\tilde{X}$ of $X$ where the inverse images of $E_1$ and $E_2$ are ruled surfaces over each elliptic curve.  The following modularity result was shown:

\begin{thm}
\cite{bib:Lee}
There exists a resolution $\tilde{X}$ of $X$ such that, up to semisimplification,

\begin{equation*}
H^3_{et}(\tilde{X}, \QQ_l) \cong U \oplus \Ind_{Gal(\overline{\QQ}/\QQ(i))}^{Gal(\overline{\QQ}/{\QQ})} H^1((E_1), \QQ_l)(-1)
\end{equation*}

where $U$ is a two-dimensional Galois representation unramified outside the bad primes $\{2,5 \}$ of $X$.  The representation $U$ corresponds to the unique normalized cusp form $f_5$ of weight 4 and level 5.
\end{thm}

The present example is notable because of the presence of elliptic surfaces which are not defined over $\QQ$ (only their scheme-theoretic union is) and also because the level of the modular form does not divide all of the bad primes, as was conjecturally the case in the previous example of Sch\"utt.

\section{Kummer-type constructions}

Recently several examples of modular non-rigid Calabi-Yau threefolds have been obtained via Kummer-type constructions.

\subsection{Quotients of type $(K3 \times E)/<-1,-1>$}

In \cite{bib:LY} Livn\'e and Yui construct four classes of non-rigid Kummer-fibred Calabi-Yau threefolds which reach the strong Arakelov-Yau bound.  They start by considering an elliptic modular surface $Y$ associated to a torsion-free congruence subgroup $\Gamma$ of $PSL_2(\ZZ)$ of genus 0 and index 24 and lifting them to subgroups $\tilde{\Gamma}$ of $SL_2(\ZZ)$.  Let $Y(\tilde{\Gamma})$ be the elliptic modular surface associated to $\tilde{\Gamma}$; it is an elliptic fibration over the genus-zero modular curve $X(\tilde{\Gamma})$.  Now let $E$ be a elliptic curve defined over $\QQ$, and consider the threefold $Y(\tilde{\Gamma}) \times E$ as a fibration of abelian surfaces over $X(\tilde{\Gamma})$.  Let $X(\tilde{\Gamma}) \rightarrow X(\tilde{\Gamma})$ be the associated Kummer family obtained by quotienting out by the involution $-1$ on each fiber and blowing up the fixed locus.  Livn\'e and Yui prove that

\begin{thm}
In the cases $\Gamma = \Gamma(4), \Gamma_0(3) \cap \Gamma(2), \Gamma_0(8) \cap \Gamma(2), \Gamma_1(8;4,1,2)$, the above construction yields a non-rigid Calabi-Yau threefold $X(\tilde{\Gamma})$ whose generic fiber is a Kummer surface and which achieves the strong Arakelov-Yau bound

\begin{equation*}
\deg f_{*} \omega_{X/\PP^1} \leq \frac{1}{2} \deg \Omega^1_{\PP^1} (\log S),
\end{equation*}

where $S$ is the set of cusps of $\Gamma$.  The $L$-series of $X$ is

\begin{equation*}
\begin{aligned}
L(H^3(X),s) &= L(T(Y) \times H^1(E),s) L(N_{-} \times H^1(E), s) \\
&= L(g_3 \otimes g_2, s) L(g_E, s-1)^k \prod_{\delta} L(g_E \otimes \chi_{\delta}, s-1)^{k_{\delta}}. \\
\end{aligned}
\end{equation*}
\end{thm}

Here $N_{-}$ is the motive of algebraic cycles on $Y$ which are odd under the action of the involution $-1$, $g_3$ is the weight 3 modular form associated to the rank 2 motive $T$, $g_E$ is the weight 2 modular form of the elliptic curve $E/ \QQ$, $k$ is the number of cycles in $N_{-}$ defined over $\QQ$, $k_{\delta}$ is the number of cycles in $N_{-}$ defined over $\QQ(\sqrt{\delta})$ and $\chi_{\delta}$ is the quadratic character cut out by $\QQ(\sqrt{\delta})$.

The proof depends on the fact that for the given groups $\Gamma$, $Y(\tilde{\Gamma})$ is an singular $K3$ surface:  it has Picard number 20, as discussed in Section 3.  Define the rank $2$ motive of transcendental cycles on $Y(\tilde{\Gamma})$ to be $T = H^2(Y(\tilde{\Gamma}))/ NS(Y(\tilde{\Gamma}))$.  Then the modular form $g_3$ of weight 3 is the modular form associated to the rank $2$ motive $T$.  The intuition here is that $H^3(X)$ is the part of $H^3(Y \times E)$ which is even under the involution $(-1,-1)$, and this is equal to (odd part of $H^2(Y)$ $\otimes$ odd part of $H^1(E)$ $\oplus$ even part of $H^2(Y)$ $\otimes$ even part of $H^1(E)$), since $H^1(Y) = 0$.  Since the holomorphic 1-form on $E$ is odd, $H^1(E)$ has no even part, and we are left with (odd part of $H^2(Y)$ $\otimes$ odd part of $H^1(E)$).  The odd part of $H^2(Y)$ consists of the transcendental part of $Y$ (the holomorphic 2-form) plus the odd algebraic cycles on $Y$.

Borcea and Voisin \cite{bib:Borcea}, \cite{bib:Voisin} independently gave a construction of mirror families of Calabi-Yau threefolds of type $K3 \times E / <-1,-1>$.  It would be interesting to compare the $L$-series and zeta functions of modular threefolds in mirror families and compare them with recent results of Wan \cite{bib:Wan}.

\subsection{A quotient of the Schoen quintic}

In this section we consider a threefold $X$ obtained by resolving a $\ZZ/2$ quotient of the Schoen quintic threefold $Q$ in $\PP^4$ given by the equation $x_0^5 + x_1^5 + x_2^5 + x_3^5 + x_4^5 - 5x_0x_1x_2x_3x_4 = 0$; it does not seem to have been studied in the current literature.  In \cite{bib:Schoen1} it was shown that $Q$ has the following properties:

\begin{enumerate}

\item $Q$ is the total space of a $\PP^1$ family of abelian surfaces arising from the Horrocks-Mumford vector bundle on $\PP^4$;

\item $Q$ is modular in the sense that the middle-dimensional cohomology group of a resolution of $Q$ has $L$-series equal to the Mellin transform of a weight 4 modular form of level 25, up to an Euler factor at the bad prime 5.

\end{enumerate}

We show that $X$ has the following properties:

\begin{enumerate}

\item $X$ is a Calabi-Yau variety;

\item $X$ is the total space of a $\PP^1$ family of Kummer surfaces;

\item A big resolution $\tilde{X}$ of $X$ is modular, and the middle-dimensional cohomology group of $\tilde{X}$ is, up to semisimplification, the direct sum of $H^3(Q)$ and $H^3(S)$ where $S$ is a $\PP^1$-bundle over an elliptic curve $E$ defined over $\QQ$.

\end{enumerate}

Our starting point is the observation in \cite{bib:HM} that the zero set of a generic section $s$ of $HM$ is an abelian surface $A_s$ of degree 10, and that the involution $\iota$ on $\PP^4$ defined by $\iota[x_0:x_1:x_2:x_3:x_4] = [x_0:x_4:x_3:x_2:x_1]$ induces the negation morphism on $A_s$ if the origin of $A_s$ is chosen to be one of the fixed points of $\iota$.

The involution $\iota$ has fixed loci $\PP^1_{-} = \{ x_0 = x_1+x_4 = x_2+x_3 = 0 \}$ and $\PP^2_{+} = \{ x_1 - x_4 = x_2 - x_3 \} = 0$.  It was proven by Horrocks-Mumford that the restriction map $\Gamma(HM) \rightarrow \Gamma(HM_{\PP^1_{-}})$ is an injection, that $\Gamma(HM_{\PP^1_{-}}) \cong \OO(-1) \oplus \OO(6)$ and that the image of this restriction map is spanned by the sections $\lambda(\lambda^5 + 2 \mu^5), \mu(\mu^5 - 2\lambda^5), 5\lambda^4\mu^2, 5\lambda^2\mu^4$.  Hence the generic $A_s$ intersects $\PP^1_{-}$ in 6 points.  $A_s$ intersects $\PP^2_{+}$ in 10 points, and if one of these 16 points is chosen to be the origin for $A_s$ then these 16 points are exactly the points of order 2.  Taking the quotient of $A_s$ by $\iota$ and blowing up the 16 nodes yields a Kummer surface.

A Horrocks-Mumford quintic is a threefold $Q_{s_1 \wedge s_2}$ given by a determinantal section $s_1 \wedge s_2$ of the bundle $\wedge^2 HM \cong \OO(5)$ on $\PP^4$.  One can easily check that in the generic case, the singularities of this threefold are double points at the 100 points of intersection of $A_{s_1}$ and $A_{s_2}$.  A small resolution $\hat{Q}_{s_1 \wedge s_2}$ can be obtained by considering the section of the bundle $HM \boxtimes \OO(1) \rightarrow \PP^4 \times \PP^1$ given by $c_1 s_1 + c_2 s_2$, where $[c_1 : c_2]$ are coordinates on $\PP^1$.  This small resolution of $Q_{s_1 \wedge s_2}$ is a fibration over $\PP^1$ whose generic fiber is an abelian surface.

In \cite{bib:Schoen1} it was shown that the quintic $Q$ defined by the equation $x_0^5 + x_1^5 + x_2^5 + x_3^5 + x_4^5 - 5x_0x_1x_2x_3x_4$ is a Horrocks-Mumford quintic.  It is a special member of the family of HM-quintics, having 125 nodes at the points $[1:\zeta^a:\zeta^b:\zeta^c:\zeta^{-c-b-a}]$ where $\zeta$ is a primitive 5th root of unity, instead of just the usual 100.

The fixed locus of the involution $\iota$ on $Q$ consists of the line $\PP^1_{-}$ and the curve $E = \PP^2_{+} \cap Q = \{ [x_0:x_1:x_2] : x_0^5 + 2 x_1^5 + 2 x_2^5 - 5 x_0 x_1^2 x_2^2 = 0 \}$.  The images of the five nodes $[1: \zeta^{a} : \zeta^{-a}: \zeta^{-a} : \zeta^a]$ lie on $E$.  Since $E$ is a curve of arithmetic genus 6 with 5 nodes, its normalization is an elliptic curve.  Of course, $Q / \iota$ still has $\frac{125-5}{2} = 60$ nodes away from $\PP^1_{-}$ and $E$.

For computational purposes, it is convenient to consider a linear change of coordinates $[y_0:y_1:y_2:y_3:y_4] = [x_0:x_1 + x_4: x_2+x_3: x_2-x_3: x_1 - x_4]$ on $\PP^4$.  The Schoen quintic then has the equation 

\begin{equation*}
16 y_0^5 + (y_1^5 + 10y_1^3 y_4^2 + 5 y_1 y_4^4) + (y_2^5 + 10 y_2^3 y_3^2 + 5 y_2 y_3^4) - 5 y_0 (y_1^2 - y_4^2)(y_2^2 - y_3^2) = 0.
\end{equation*}

In $y$-coordinates, the action of $\iota$ is given by the matrix

\begin{equation*}
\begin{pmatrix} 1 & & & & \\ & 1 & & & \\ & & 1 & & \\ & & & -1 & \\ & & & & -1 \\ \end{pmatrix}
\end{equation*}

The homogeneous ring of invariants of $\PP^4$ is simply $k[y_0, y_1, y_2, y_3^2, y_4^2, y_3 y_4]$ = $k[Y_0, Y_1, Y_2, Y_3, Y_4, Y_5] / Y_3 Y_4 - Y_5^2$.  Hence $\PP^4/\iota$ can be exhibited as the hypersurface $Y_3 Y_4 - Y_5^2 = 0$ in the toric variety $\PP(1,1,1,2,2,2)$, and $Q/\iota$ is a complete intersection in $\PP(1,1,1,2,2,2)$ given by the equations

\begin{equation*}
\begin{split}
16 Y_0^5 + (Y_1^5 + 10 Y_1^3 Y_4 + 5 Y_1 Y_4^2) + (Y_2^5 &+ 10 Y_2^3 Y_3 + 5 Y_2 Y_3^2) \\
&- 5 Y_0 (Y_1^2 - Y_4)(Y_2^2 - Y)3) = 0 \\
&Y_3 Y_4 - Y_5^2 = 0 \\
\end{split}
\end{equation*}

Simple but messy calculations show that the singularities of $Q/\iota$ along $\PP^1_{-}$ are transverse $A_1$ singularities, and that the singularities along $E$ are also transverse $A_1$ singularities except at the 5 special points of $E$, which are $A_3$ singularities (blowing up $E$ yields a $\PP^1$-bundle over $E$ which still has 10 double points lying over the original 5 singularities on $E$).  We blow up the $\PP^1_{-}$ and blow up the ambient $\PP^2_{+}$ containing $E$; the resulting threefold $X$ has the original 60 nodes from $Q/\iota$ and 10 more nodes lying in the exceptional divisor over $E$. After making an appropriate choice of small resolution of the 60+10 nodes, we obtain a smooth threefold $X$ which is a fibration over $\PP^1$ with generic fiber a Kummer surface.

Working over $\CC$, the holomorphic 3-form on $Q$ (away from the nodes) is given by 

\begin{equation*}
\Omega = \Res \frac{\sum_j (-1)^j x_j dx_0 \wedge \dots \wedge \hat{dx_j} \wedge \dots \wedge dx_n}{x_0^5 + \dots + x_4^5 - 5x_0x_1x_2x_3x_4}.
\end{equation*}

It is simple to check that $\Omega$ is invariant under $\iota$ and thus descends to a nonvanishing holomorphic 3-form on $Q/\iota$.  Since nonvanishing holomorphic 3-forms extend across transverse $A_1$-singularities and small resolutions of nodes, we obtain a nonvanishing 3-form on $X$.  Hence $X$ is Calabi-Yau.  The bad primes of $X$ are checked to be 2 and 5.

Now let $\tilde{X}$ be a big resolution of $X$, where we replace the exceptional $\PP^1$s over nodes with $\PP^1 \times \PP^1$.  We now compute the Euler characteristic of $\tilde{X}$. A smooth quintic threefold has Euler number -200.  $Q$ is obtained by shrinking 125 $S^3$ cycles to points, so the Euler number of $Q$ is -75.

$Q$ is a branched cover of $Q/\iota$ which is generically 2:1 and branched over $\PP^1_{-} \cup E$.  Now $E$ is a curve of arithmetic genus 6 with 5 $S^1$-cycles shrunk to points, so its Euler number is -5, while the Euler number of $\PP^1_{-}$ is 2.  By the Riemann-Hurwitz formula, $2 \chi(Q/\iota) - (2 - 5) = -75$, so the Euler characteristic of $Q/\iota$ is -39.

$\tilde{X}$ is obtained from $Q/\iota$ by replacing $\PP^1_{-}$ with a $\PP^1$-bundle $S$ over $\PP^1_{-}$, replacing $E$ with a $\PP^1$-bundle over $E$, and taking a big resolution of a total of $60 + 5 \cdot 2 = 70$ nodes.  The resulting Euler characteristic of $\tilde{X}$ is 168.

To compute the Hodge numbers of $\tilde{X}$, we use a trick due to van Geemen and Werner \cite{bib:vGW} (which has in fact been implicitly used in some of the examples discussed above).  Applying the Lefschetz theorem to the Frobenius morphism, we have

\begin{equation*}
\begin{aligned}
\#X_p &= \sum_{i=0}^6 (-1)^i \trace \Frob_p(H^i(\tilde{X}))
\end{aligned}
\end{equation*}

Here $\# \tilde{X}_p$ is the number of fixed points of $X$ under the Frobenius morphism, i.e. the number of $\FF_p$-rational points of $X$.

For $p \equiv 1 ({\rm mod} 20)$, $H^2(\tilde{X})$ is spanned by algebraic cycles defined over $\FF_p$:  the only thing to check is that the rulings on the exceptional $\PP^1 \times \PP^1$s over the nodes are defined over $\FF_p$.  The coordinates of the nodes are defined whenever $i$ is in $\FF_p$, and local calculations show that the rulings are defined whenever $-5$ is a square in $\FF_p$.  Hence  $\trace \Frob_p (H^2(X)) = p h^2$, where $h^2$ is the second Betti number of $X$.  By Poincar\'e duality, $\trace \Frob_p (H^4(X)) = p^2 h^2$.  Now we know that $|\trace \Frob_p (H^3(X))| \leq p^{\frac{3}{2}} h^3 = p^{\frac{3}{2}} (168 - 2 - 2h^2)$.  By counting points (using a computer) on $\FF_p$, we can hope to determine $h^2$.  In fact, choosing $p = 421$ yields $h^2(X) = 85$ and $h^3(X) = 4$.

Consider again the curve $E = \{ x_0^5 + 2x_1^5 + 2x_2^5 - 5x_0 x_1^2 x_2^2 = 0 \}$ in $\PP^2_{+}$.  It has nodes at $[1: \zeta^a : \zeta^{-a}]$.  One can count points on the normalization $\hat{E}$ of $E$ over $\FF_p$:  one simply needs to keep track of whether the nodes and their branches are defined over $\FF_p$.  As before, counting points on $E$ enables us to determine the values of $a_p = \trace \Frob_p (H^1(\hat{E}))$; the bad primes of $\hat{E}$ are 2 and 5.  

We want to find the modular form associated to $\hat{E}$.  For good primes $p < 100$, the values of $a_p$ are the same as the values of the newform $50B1$ of level 50 and weight 2.  We would like to say that the Galois representations associated to $E$ and to the newform $50B1$ are isomorphic, but the traces are not even.  However, by the method of Serre explained above and worked out in \cite{bib:Schoen1}, a two-dimensional Galois representation ramified at $\{2,5 \}$ is uniquely determined by its determinant and its traces at $p \in T = \{3,7,11,13,17,29,31 \}$.  The traces do indeed agree at this set of primes $T$, so we have found the correct modular form.

We have a rational map from $\PP^1$-bundle $S$ over $\hat{E}$ into $\hat{X}$ which can be resolved by blowups.  Hence it induces a nontrivial map on $H^3$, so we have a decomposition (up to semisimplification)

\begin{equation*}
\begin{aligned}
H^3(X) &\cong H^3(S) \oplus L \\
&\cong H^1(\hat{E})(-1) \oplus L.
\end{aligned}
\end{equation*}

Now the quotient map defines a correspondence between $Q$ and $\hat{X}$.  Since the holomorphic three-form on $Q$ is preserved by the involution, the quotient operation preserves the motive $H^3(Q)$.  Hence the same level 25 newform of $Q$ must occur in the decomposition of $H^3(\tilde{X})$.  We therefore have

\begin{thm}
The threefold $\tilde{X}$ is modular, and up to Euler factors at the bad primes 2 and 5, the $L$-series of $X$ is given by

\begin{equation}
L(H^3(\hat{X},s)) = L(f_{25},s) L(g_{50}, s-1)
\end{equation}

where $f_{25}$ is a weight 4, level 25 newform and $g_{50}$ is a weight 2, level 50 form.

\end{thm}

\section{Double octics}

In \cite{bib:CM} Cynk and Meyer considered examples of Calabi-Yau threefolds constructed as double covers of $\PP^3$ branched over a degenerate octic (an arrangement of eight planes) $D$ and proved their modularity by several of the methods used above; further computation was done in \cite{bib:Meyer}.  Using a computer search, they found degenerate arrangements of $D$ with $h^{12} = 1$ and found candidate modular forms for the resulting threefolds.  Some were proven to be modular by finding an elliptic surface $S$ in $X$ and thus splitting $H^3(X)$ into an elliptic piece and a piece $U$ corresponding to a weight 4 modular form, using the technique of Hulek-Verrill.  Others were proven to be modular by exhibiting birational maps between them and known modular rigid Calabi-Yaus or the Kummer-type non-rigid Calabi-Yaus of Livn\'e-Yui.

For some examples Cynk and Meyer found an involution of $\PP^3$ preserving the arrangement $D$, leading to an involution on the resulting Calabi-Yau threefold $X$.  They then use a trick of van Geemen and Nygaard \cite{bib:NvG} to split $H^3$ into two two-dimensional pieces:  if both 1 and -1 are eigenvalues of the action of the involution on $H^3$, then $H^3$ splits.  Since the splitting commutes with the Frobenius morphism, it splits the Galois representation $H^3$ into two two-dimensional subrepresentations.

If $\Phi$ is a $\QQ$-linear involution on $\PP^3$ which preserves the arrangement of planes $D$, then it induces an involution on the resulting Calabi-Yau $X$.  One can apply the Lefschetz theorem to $\Frob_p \circ \Phi$.  The fixed points of $\Frob_p \circ \Phi$ are simply the $\FF_{p^2}$-rational points of $X$, which can be counted via computer.  The map $\Frob_p \circ \Phi$ acts by multiplication by $\pm p$ on $H^2$ and by $\pm p^2$ on $H^4$.  The trace of $\Frob_p \circ \Phi$ will be $\pm 1$ times the trace of $\Frob_p$ on the $\pm 1$ eigenspaces of $H^3$.  Thus computing the traces of $\Frob_p \circ \Phi$ by counting points on $\FF_{p^2}$ and computing the traces of $\Frob_p$ will yield the traces of $\Frob_p$ on the $\pm 1$ eigenspaces of $H^3$.  Cynk and Meyer show that for some of their examples, the group $H^3$ does indeed split via the action of some involution $\Phi$, and the two-dimensional pieces are proven to be modular using the Livn\'e theorem.

\subsection{Higher-dimensional Kummer constructions}

Higher-dimensional modular Calabi-Yau varieties have been considered by Cynk-Hulek and Cynk-Sch\"utt. 

In \cite{bib:CH} Cynk-Hulek consider products of elliptic curves with automorphisms of finite order:

\begin{thm}
Let $E_1, E_2, \dots, E_n$ be elliptic curves defined over $\QQ$ with involutions on each given as $x \rightarrow -x$.  Consider the subgroup $\ZZ_2^{n-1} \subset \ZZ_2^n$ of vectors $(c_1, c_2, \dots, c_n)$ with $\sum c_i = 0$, and put $X^n =  (E_1 \times \dots \times E_n)/\ZZ_2^{n-1}$.  Then $X^n$ admits a crepant resolution $\hat{X}^n$ which is a smooth Calabi-Yau variety.

If $n$ is odd, then $H^n(\hat{X}^n) \cong H^n(E_1 \times \dots \times E_n)^{\ZZ_2^{n-1}} \cong H^1(E_1) \otimes \dots \otimes H^1(E_n)$.  Then $\hat{X}^n$ is modular, with $L(H^n(\hat{X}^n), s) = L(g_{E_1} \otimes \dots \otimes g_{E_n}, s)$, where the $g_{E_i}$ are the modular forms associated to the $E_i$.

If $n$ is even, then the invariant submotive $H^n(E_1 \times \dots \times E_n)^{\ZZ_2^{n-1}}$ is isomorphic to the direct sum of a submotive generated by cycles of products of $n/2$ fibers and a submotive $I(\hat{X}^n) \cong H^1(E_1) \otimes \dots \otimes H^1(E_n)$, and $I(\hat{X}^n)$ contains the submotive of transcendental cycles.  In this case $L(I(\hat{X}^n),s) = L(g_{E_1} \otimes \dots \otimes g_{E_n},s)$.
\end{thm}

If one wants a rigid Calabi-Yau variety, one must divide out by automorphisms of higher order:

\begin{thm}
Let $E$ be the elliptic curve with automorphism $\eta$ of order 3, and let $X^n$ be the quotient of $E^n$ by the group $\{ (\eta^{a_1},\dots,\eta^{a_n}) : \sum a_i = 0 \}$.  Then $X^n$ has a crepant resolution $\hat{X}^n$ which is Calabi-Yau, and $\dim H^n(\hat{X}^n) = 2$ if $n$ is odd, resp. $\dim T(\hat{X}^n) = 2$ if $n$ is even, where $T(\hat{X}^n)$ is the transcendental part of the middle cohomology.

Moreover, $\hat{X}^n$ is defined over $\QQ$ and $L(H^n(\hat{X}^n), s) = L(g_{n+1},s)$ resp. $L(T(\hat{X}^n),s) = L(g_{n+1},s)$ where $g_{n+1}$ is the weight $(n+1)$ cusp form with complex multiplication in $\QQ(\sqrt{-3})$ associated to the $n$th power of the Grossencharacter of $E$.
\end{thm}

\begin{thm}
Let $E$ be the elliptic curve with Weierstrass equation $y^2 = x^3 - Dx$, where $D$ is a square-free integer, $E$ has complex multiplication in $\QQ(i)$ and the map $\rho: (x,y) \rightarrow (-x,iy)$ is an automorphism of order 4.

Let $X^n$ be the quotient of $E^n$ by the group $\{ (\eta^{a_1} \times \dots \times \eta^{a_n}) : \sum a_i = 0 \}$.  Then $X^n$ has a crepant resolution $\hat{X}^n$ defined over $\QQ$ which is Calabi-Yau.  We have $\dim H^n(\hat{X}^n) = 2$ if $n$ is odd, resp. $\dim T(\hat{X}^n) = 2$ if $n$ is even, where $T(\hat{X}^n)$ is the transcendental part of the middle cohomology.

Then $L(H^n(\hat{X}^n),s) = L(g_{n+1},s)$ resp. $L(T(\hat{X}^n),s) = L(g_{n+1},s)$ where $g_{n+1}$ is a weight $n+1$ cusp form with complex multiplication in $\QQ(i)$.
\end{thm}

In \cite{bib:CS} Cynk and Sch\"utt generalize this higher-dimensional Kummer construction.  Let $L$ be a number field, $L(\alpha)/L$ an extension of order $n$, $E_1$ an elliptic curve over $L(\alpha)$, and $E_2, \dots, E_n$ the conjugates of $E_1$.  Construct the product $A = E_1 \times \dots \times E_n$ and the Kummer quotient $X = A/ \ZZ_2^{n-1}$ as before.  They then show that $A$ and $X$ descend to the base field $L$ (Weil restriction).  In the case $n = 3$, they construct a crepant resolution $\tilde{X}$ of $X$ defined over the base field $L$.

Now suppose $K$ is a quadratic imaginary field of class number 3, $d = -p_0$ its discriminant, $E_1$ an elliptic curve with complex multiplication in $\OO_K$, $F$ the field of definition of $E$.  Construct the Calabi-Yau threefold $\tilde{X}$ as above.  We then have

\begin{thm}
The Calabi-Yau threefold $\tilde{X}$ is modular, and

\begin{equation*}
L(\tilde{X}/\QQ,s) = L(g,s) L(E/F,s-1)
\end{equation*}

where $g$ is the unique $CM$-newform of level $d^2$ and weight 4 with rational coefficients, or equivalently the corresponding Hecke character $\phi$ of $K$ with conductor $(\sqrt{d})$ and $\infty$-type 3.

\end{thm}

\section{A Hilbert-modular Calabi-Yau threefold}

In \cite{bib:CoS} Consani and Scholten study a quintic threefold $X$ given by the homogenization of the equation $P_5(x,y) - P_5(z,w)$ where

\begin{equation*}
P_5(x,y) = (x^5+y^5) - 5xy(x^2+y^2) + 5xy(x+y) + 5(x^2+y^2) - 5(x+y)
\end{equation*}

This Calabi-Yau has 120 ordinary double points as singularities.  Let $\tilde{X}$ be the big resolution of $X$.  It is shown that $h^3(\tilde{X}) = 4$ and $\tilde{X}$ has bad reduction at $\{2,3,5\}$.  They prove that the Galois representation $\rho$ associated to $H^3_{et}(\tilde{X}, \QQ_l(\sqrt{5}))$ splits into two two-dimensional pieces.  They also construct a Hilbert modular newform $f$ of weight $(2,4)$ and conductor 30 on a finite extension $E_\lambda$ of $\QQ_l$ with ring of integers $\OO_{\lambda}$ and an associated 2-dimensional $\lambda$-adic Galois representation $\sigma_{f,\lambda} \rightarrow GL_2(\OO_{\lambda})$ and give numerical evidence that $H^3_{et}(\tilde{X}, \QQ_l(\sqrt{5}))$ is associated to the forms $f$ and $f^{\sigma}$.  In 2005 Yi announced a proof of this conjecture \cite{bib:Yi}.

\section{Further questions}

We conclude this paper with some open problems.

\begin{question}
For which weight 4 newforms $f$ does there exist a Calabi-Yau threefold $X$ such that $L(f,s)$ occurs in the $L$-series of $X$? 
\end{question}

This question was posed independently by Mazur and van Straten.  When $f$ is a newform of weight 2 with coefficents in the field $F$, a construction of Shimura \cite{bib:Shimura} associates to $f$ an abelian variety $A/\QQ$ of dimension $[F:\QQ]$.  No such construction yet exists for higher weight.  The known examples typically arise as highly symmetric members in families of Calabi-Yau varieties.

\begin{question}
If $X$ is a modular Calabi-Yau threefold, can one predict the level of the associated modular form $f$?
\end{question}

In the elliptic curve case, the level of the newform is given by the conductor $N$, which is described in terms of the bad reduction of $N$.  However, we have seen examples of Calabi-Yau threefolds where the bad primes do not appear in the level of the modular form.  Meyer has formulated the following rule of thumb:

\begin{conj}  
Let $X$ be a modular Calabi-Yau threefold and let $p \geq 5$ be a prime.  Let $f$ be the twist of minimal level $N$ of the weight four newform associated with $X$.  If the exponent of $p$ in $N$ is 2 then $X$ develops non-isolated singularities over $\FF_p$.  If the exponent is 1 then $X$ is singular modulo $p$ and there exists a birational model of $X$ with only isolated singularities modulo $p$.  If the exponent is 0 then there is a birational model of $X$ with good reduction at $p$.
\end{conj}

Hulek, Kloosterman and Sch\"utt also pose the question

\begin{question}
\cite{bib:HKS}
Let $X$ be a modular Calabi-Yau variety over $\QQ$ with corresponding newform $f$.  Under what condition do the bad Euler factors of the $L$-series of $X$ and $f$ coincide?
\end{question}

\begin{question}
More generally, what other kinds of automorphic forms can occur in the $L$-series of a Calabi-Yau threefold over $\QQ$?  What about over other number fields?
\end{question}


\end{document}